\documentstyle[12pt]{article}

\input amssym.def
\input amssym

\def\twelveamsfonts{
 \font\twelvemsa=msam10 scaled 1200
 \font\eightmsa=msam8
 \font\sixmsa=msam6
 \font\twelvemsb=msbm10 scaled 1200
 \font\eightmsb=msbm8
 \font\sixmsb=msbm6
 \font\twelvembi=cmmib10 scaled 1200
 \font\eightmbi=cmmib8
 \font\sixmbi=cmmib6
 \textfont\msafam\twelvemsa
 \scriptfont\msafam\eightmsa
 \scriptscriptfont\msafam\sixmsa
 \textfont\msbfam\twelvemsb
 \scriptfont\msbfam\eightmsb
 \scriptscriptfont\msbfam\sixmsb}

\ifnum\ht0<\ht2\let\twelveamsfonts\relax\fi
\twelveamsfonts

\let\bls\baselineskip \let\nt\noindent
\let\vp\vphantom 
  
\def\vsk#1>{\vskip#1\bls} \def\vv#1>{\vadjust{\vsk#1>}}
\def\,{\relax\ifmmode\mskip\thinmuskip\relax\else\kern.16667em\fi}
\def\;{\relax\ifmmode\mskip\thickmuskip\relax\else\kern.27777em\fi}
\def\!{\relax\ifmmode\mskip-\thinmuskip\relax\else\kern-.16667em\fi}
\def\&{.\kern.1em} \def\itl#1{{\it #1\/}} 
\def\ftext#1{{\let\thefootnote\relax\footnotetext{\vsk-.8>\noindent #1}}}

\let\ge\geqslant
\let\geq\geqslant
\let\le\leqslant
\let\leq\leqslant

\def\bea{\begin{eqnarray}}
\def\ena{\end{eqnarray}}
\def\beq{\begin{equation}}
\def\eeq{\end{equation}}
\def\no{\nonumber}
\def\qed{\quad$\square$}

\def\proof{\noindent{\it Proof.}\quad}

\newtheorem{thm}{Theorem}[section]
\newtheorem{prop}[thm]{Proposition}

\makeatletter
\newbox\p@b@ld
\def\poorbold#1{\setbox\p@b@ld\hbox{#1}\kern-.01em\copy\p@b@ld\kern-\wd\p@b@ld
 \kern.02em\copy\p@b@ld\kern-\wd\p@b@ld\kern-.012em\raise.02em\box\p@b@ld}

\@addtoreset{thm}{section}
\@addtoreset{equation}{section}
\makeatother

\begin{document}

\begin{center}
\vp1
{\Large \bf On solutions of the $q$-hypergeometric equation
\vsk.3>
with $q^{N}=1$
}
\vsk2>
{Yoshihiro Takeyama$^{\,\diamond}$}
\ftext{$^{\diamond\,}$Research Fellow of the Japan Society for
the Promotion of Science. \\
e-mail:ninihuni@kurims.kyoto-u.ac.jp
}
\vsk1.5>
{\it Research Institute for Mathematical Sciences,
Kyoto University, Kyoto 6068502, Japan}
\end{center}
\vsk1.75>

{\narrower\nt
{\bf Abstract.}\enspace
We consider the $q$-hypergeometric equation with $q^{N}=1$ and 
$\alpha, \beta, \gamma \in {\Bbb Z}$. 
We solve this equation on the space of functions 
given by a power series multiplied by a power of the logarithmic function. 
We prove that the subspace of solutions is two-dimensional 
over the field of quasi-constants. 
We get a basis for this space explicitly. 
In terms of this basis, 
we represent the $q$-hypergeometric function of the Barnes type 
constructed by Nishizawa and Ueno. 
Then we see that this function has logarithmic singularity at the origin. 
This is a difference between the $q$-hypergeometric functions 
with $0<|q|<1$ and at $|q|=1$.
\vsk1.4>}
\vsk0>
\thispagestyle{empty}

\section{Introduction}
Consider the $q$-hypergeometric equation
\bea
\left\{
(1-D_{q})(1-q^{\gamma-1}D_{q})-t(1-q^{\alpha}D_{q})(1-q^{\beta}D_{q})
\right\}
\varphi(t)=0,
\label{uqhyp}
\ena
where $D_{q}$ is the $q$-difference operator 
defined by $(D_{q}\varphi)(t):=\varphi(qt)$. 
In this paper, we solve (\ref{uqhyp}) with $q^{N}=1$ and $\alpha, \beta, \gamma \in {\Bbb Z}$ 
explicitly on a special space of functions 
and represent the $q$-hypergeometric function of the Barnes type at $|q|=1$ 
in terms of our solutions. 

Let us recall some results about solutions to (\ref{uqhyp}). 
In the case of $0<|q|<1$, one of the solutions to (\ref{uqhyp}) is 
the basic hypergeometric function $\varphi(\alpha, \beta, \gamma ; t)$ \cite{GR} 
defined by
\bea
\varphi(\alpha, \beta, \gamma ; t):=
\sum_{k=0}^{\infty}
\frac{(q^{\alpha})_{k}(q^{\beta})_{k}}{(q)_{k}(q^{\gamma})_{k}}t^{k},
\label{defbhyp}
\ena
where $(x)_{n}:=\prod_{j=0}^{n-1}(1-q^{j}x).$ 
We can get this solution by setting
\bea
\varphi(t)=\sum_{k=0}^{\infty}c_{k}t^{k}, \quad c_{k} \in {\Bbb C}
\label{powerseries}
\ena
and solving a recursion relation for $\left\{ c_{k} \right\}$. 
In the case of $|q|=1$, we can not get solutions in this manner 
because the coefficient in (\ref{defbhyp}) does not converge. 
However, some solutions are constructed in terms of a contour integral 
by Nishizawa and Ueno \cite{NU}. 
In this paper, we compare the solutions of these two types 
under the condition that $q^{N}=1$ and $\alpha, \beta, \gamma \in {\Bbb Z}$. 

First we try to find solutions in a similar way to the case of $0<|q|<1$. 
Then we consider solutions of more general form than (\ref{powerseries}) 
because of the following reason.
In the limit as $q \to 1$, the equation (\ref{uqhyp}) goes to 
the hypergeometric differential equation:
\bea
\left\{
t(1-t)\frac{d^{2}}{dt^{2}}+(\gamma-(\alpha+\beta+1)t)\frac{d}{dt}-\alpha\beta
\right\}F(t)=0.
\ena
It is known that if $\gamma \not\in {\Bbb Z}$ 
there exist two independent solutions at $t=0$ of the form
\bea
F(t)=\sum_{k=0}^{\infty}a_{k}t^{k+\rho}.
\label{class}
\ena
However, if $\gamma \in {\Bbb Z}$ one of the two independent solutions 
is represented as
\bea
F_{1}(t)\log{t}+F_{2}(t)
\label{logclass}
\ena
where $F_{1}(t)$ and $F_{2}(t)$ are functions of the from (\ref{class})
(see \cite{AAR}, for example).

Now let us return to the equation (\ref{uqhyp}) with $0<|q|<1$. 
When we consider solutions of the form (\ref{class}), 
we get two independent solutions if $\gamma \not \in {\Bbb Z}$. 
One of them is the basic hypergeometric function (\ref{defbhyp}) 
and the other is given by
\bea
t^{1-\gamma}\varphi(\alpha-\gamma+1, \beta-\gamma+1, 2-\gamma ; t)
=\sum_{k=0}^{\infty}
\frac{(q^{\alpha-\gamma+1})_{k}(q^{\beta\gamma+1})_{k}}{(q)_{k}(q^{2-\gamma})_{k}}
t^{k+1-\gamma}.
\ena
Here we note that if $\gamma \in {\Bbb Z}$ one of these solutions does not make sense. 
Then we can construct another solution in the form (\ref{logclass}) 
(see \cite{E} for details).

{}From the consideration above, 
we formulate our problem as follows.
We set $t=q^{x}$ and 
rewrite (\ref{uqhyp}) as a difference equation for a function of $x$, see (\ref{qhyp}). 
We try to find solutions of the following form:
\bea
\sum_{j=0}^{n}\sum_{k=0}^{\infty}c_{jk}x^{j}q^{kx}, 
\quad c_{jk} \in {\Bbb C} .
\label{first}
\ena
Here the part of $j>0$ corresponds to the term in (\ref{logclass}) 
with logarithmic singularity. 
Now we note that the function $q^{Nx}$ is invariant under the shift $x \mapsto x+1$. 
Hence we can rewrite (\ref{first}) as 
\bea
\sum_{j=0}^{n}\sum_{k=0}^{N-1}f_{jk}(x)x^{j}q^{kx}, 
\label{second}
\ena
where $f_{jk}(x)$ is a periodic function with a period 1, 
that is a quasi-constant for the difference equation (\ref{qhyp}). 
We solve the $q$-hypergeometric equation 
with $q^{N}=1$ and $\alpha, \beta, \gamma \in {\Bbb Z}$ 
on the space of functions of the form (\ref{second}).

The result is as follows (Theorem \ref{th2}). 
The space of solutions is two-dimensional over the field of quasi-constants, 
and all the solutions are represented as (\ref{logclass}), 
that is $n=0$ or $1$ in (\ref{second}). 

Next we consider solutions to (\ref{uqhyp}) with $|q|=1$. 
One of the solutions is the $q$-hypergeometric function 
of the Barnes type at $|q|=1$, see (\ref{defhyp}). 
We can deal with the integral representation of this function 
in the framework of the $q$-twisted cohomology at $|q|=1$ \cite{T}. 
It is shown that, 
if $q$ is not a root of unity and the parameters $\alpha, \beta$ and $\gamma$ are generic, 
we can construct two independent solutions to (\ref{uqhyp}) 
in terms of the integral of the Barnes type 
by taking two independent homologies. 
However, in the case that $q^{N}=1$ and $\alpha, \beta, \gamma \in {\Bbb Z}$, 
the function (\ref{defhyp}) is a unique solution of this form. 
In Theorem \ref{th3}, we write down an explicit formula 
for this function in terms of the basis constructed in Theorem \ref{th2}. 
Then we see that if the parameters $\alpha, \beta$ and $\gamma $ satisfy some condition, 
the $q$-hypergeometric function at $|q|=1$ has logarithmic singularity.
On the other hand, the $q$-hypergeometric function with $0<|q|<1$ 
defined by (\ref{defbhyp}) has no logarithmic singularity. 
This is a difference between the cases of $0<|q|<1$ and $|q|=1$.

\section*{Acknoledgement}
The author thanks Masaki Kashiwara and Tetsuji Miwa for valuable remarks. 
  
\section{A basis for the space of solutions}

Let $N$ be a integer with $N \ge 2$. 
The $q$-hypergeometric difference equation is defined by 
\bea
L\Psi(x)=0, \quad 
L:=(1-D)(1-q^{\gamma -1}D)-q^{x}(1-q^{\alpha}D)(1-q^{\beta}D),
\label{qhyp}
\ena
where $D$ is the difference operator defined by $D\Psi(x):=\Psi(x+1)$.
In this paper, we consider the equation (\ref{qhyp}) with 
\bea
q:=e^{\frac{2 \pi i}{N}}, \quad \alpha, \beta, \gamma \in \{ 1, \cdots , N\}, 
\quad \beta \le \alpha.
\label{condition}
\ena
Note that the equation (\ref{qhyp}) is symmetric with respect to $\alpha$ and $\beta$, 
and hence we can assume $\beta \le \alpha$ without loss of generality.

We denote by ${\cal C}$ the field of periodic meromorphic functions of $x$ with a period 1.
This field is a space of quasi-constants for (\ref{qhyp})
in the sense that if $\Psi$ is a solution to (\ref{qhyp}) then $f\Psi$ is also a solution 
for any $f \in {\cal C}$.
 
Let us find a solution $\Psi$ of the following form:
\bea
\Psi(x)=\sum_{j=0}^{n}\sum_{k=0}^{N-1} f_{jk}x^{j}q^{kx}, 
\quad f_{jk} \in {\cal C}.
\label{assume}
\ena
It is easy to see the following.
\begin{prop}\label{unique}
The expression (\ref{assume}) is unique, that is
\bea
\sum_{j=0}^{n}\sum_{k=0}^{N-1} f_{jk}x^{j}q^{kx}=0 
\quad \Longrightarrow \quad 
\forall f_{jk}=0,
\ena
where $f_{jk} \in {\cal C}$.
\end{prop}

Let ${\cal S}$ be the space of functions of the form (\ref{assume}):
\bea
{\cal S}:=\left\{
\sum_{j=0}^{n}\sum_{k=0}^{N-1} f_{jk}x^{j}q^{kx} \, | \,
n \ge 0, f_{jk} \in {\cal C} \right\}.
\ena
Set
\bea
{\cal P}:=\left\{ \sum_{k=0}^{N-1} z_{k}q^{kx} \, | \, z_{k} \in {\cal C} \right\}.
\ena
Note that $L{\cal P} \subset {\cal P}$ and 
$L{\cal S} \subset {\cal S}$ because $q^{Nx} \in {\cal C}$. 

The following result holds.

\begin{thm} \label{th2}
The subspace of solutions to (\ref{qhyp}) in ${\cal S}$ is two-dimensional over ${\cal C}$.
A basis $\{\Psi_{1}, \Psi_{2}\}$ of this space is given as follows:

1. $ \gamma \le \beta \le \alpha$ case.
\bea
\Psi_{1}(x)&=&\sum_{k=0}^{N-\alpha}\frac{b_{0}\cdots b_{k-1}}{a_{1} \cdots a_{k}}q^{kx},  \no \\
\Psi_{2}(x)&=&x\Psi_{1}(x)+
\sum_{k=1}^{N-\alpha}\frac{b_{0}\cdots b_{k-1}}{a_{1} \cdots a_{k}}
\sum_{j=1}^{k}\left( 
\frac{1-q^{\gamma+2j-1}}{a_{j}}-\frac{1-q^{\alpha+\beta+2(j-1)}}{b_{j-1}}
\right)q^{kx} \no \\
&&
-(1-q^{\beta-\alpha})
\frac{b_{0}\cdots b_{N-\alpha-1}}{a_{1}\cdots a_{N-\alpha+1}}
\sum_{k=N-\alpha+1}^{N-\beta}
\frac{b_{N-\alpha+1}\cdots b_{k-1}}{a_{N-\alpha+2}\cdots a_{k}}q^{kx} \no \\
&&
-(1-q^{\gamma-1})\sum_{k=N-\gamma+1}^{N-1}
\frac{a_{k+1}\cdots a_{N-1}}{b_{k}\cdots b_{N-1}q^{Nx}}q^{kx}.
\label{base1}
\ena

2. $ \beta < \gamma \le \alpha$ case.
\bea
\Psi_{1}(x)=\sum_{k=0}^{N-\alpha}\frac{b_{0}\cdots b_{k-1}}{a_{1} \cdots a_{k}}q^{kx}, 
\quad \Psi_{2}(x)=\sum_{k=N-\gamma+1}^{N-\beta}
\frac{b_{N-\gamma+1}\cdots b_{k-1}}{a_{N-\gamma+2}\cdots a_{k}}q^{kx}.
\label{base2}
\ena

3. $ \beta \le \alpha < \gamma$ case.
\bea
\Psi_{1}(x)&=&\sum_{k=N-\gamma+1}^{N-\alpha}
\frac{b_{N-\gamma+1}\cdots b_{k-1}}{a_{N-\gamma+2}\cdots a_{k}}q^{kx},  \no \\
\Psi_{2}(x)&=&x\Psi_{1}(x)-
(1-q^{N-\gamma+1})\sum_{k=0}^{N-\gamma}
\frac{a_{k+1}\cdots a_{N-\gamma}}{b_{k}\cdots b_{N-\gamma}}q^{kx} \no \\
&&
+\sum_{k=N-\gamma+2}^{N-\alpha}
\frac{b_{N-\gamma+1}\cdots b_{k-1}}{a_{N-\gamma+2}\cdots a_{k}}
\sum_{j=N-\gamma+2}^{k}\left(
\frac{1-q^{\gamma+2j-1}}{a_{j}}-\frac{1-q^{\alpha+\beta+2(j-1)}}{b_{j-1}} \right)q^{kx} \no \\
&&
-(1-q^{\beta-\alpha})
\frac{b_{N-\gamma-1}\cdots b_{N-\alpha-1}}{a_{N-\gamma+2}\cdots a_{N-\alpha+1}}
\sum_{k=N-\alpha+1}^{N-\beta}
\frac{b_{N-\alpha+1}\cdots b_{k-1}}{a_{N-\alpha+2}\cdots a_{k}}q^{kx}.
\label{base3}
\ena
Here we set
\bea
a_{k}:=(1-q^{k})(1-q^{\gamma-1+k}), \quad 
b_{k}:=(1-q^{\alpha+k})(1-q^{\beta+k}).
\ena
\end{thm}

\proof
Here we prove the theorem in the first case. 
The proof for the other case is similar.

Set
\bea
\Psi=\sum_{j=0}^{n}x^{j}P_{j}, \quad P_{j}:=\sum_{k=0}^{N-1}f_{jk}q^{kx} \in {\cal P}.
\ena
{}From Proposition \ref{unique}, 
we see that $L\Psi=0$ is equivalent to the following:
\bea
LP_{j}+\sum_{t=j+1}^{n}\left( t \atop{t-j} \right)L_{t-j}P_{t}=0,
\quad (j=n, n-1, \cdots , 0), \label{cond}
\ena
where
\bea
L_{k}:=2^{k}(1-q^{x})D^{2}-((1+q^{\gamma-1})-q^{x}(q^{\alpha}+q^{\beta}))D.
\ena

It is easy to solve (\ref{cond}) for $j=n$ and $n-1$.
The solution is given by
\bea
P_{n}=f_{n}\Psi_{1}, \quad P_{n-1}=f_{n-1}\Psi_{1}+nf_{n}(\Psi_{2}-x\Psi_{1}),
\label{sol}
\ena 
where $f_{n}, f_{n-1} \in {\cal C}$.

Now we prove that $n<2$. 
If $n \ge 2$, there is a solution $P_{n-2}$ to
\bea
LP_{n-2}+(n-1)(L_{1}P_{n-1}+\frac{n}{2}L_{2}P_{n})=0,
\ena
where $P_{n-1}$ and $P_{n}$ are given in (\ref{sol}).
Especially, we have
\bea
Q_{n}:=(n-1)(L_{1}P_{n-1}+\frac{n}{2}L_{2}P_{n}) \in L{\cal P}.
\label{mustrel}
\ena

On the other hand, we can see that 
\bea
L{\cal P} \ni \sum_{k=0}^{N-1} z_{k} q^{kx}, 
\quad (z_{k} \in {\cal C})
 \quad \Longleftrightarrow \quad
\left\{
\begin{array}{l}
\displaystyle
\sum_{j=N-\beta+1}^{N-\gamma}
\frac{b_{j}\cdots b_{N-\gamma}}{a_{j}\cdots a_{N-\gamma}}z_{j}+z_{N-\gamma+1}=0, 
\quad (\gamma \not=1), \\
\displaystyle
\sum_{j=N-\beta+1}^{N-1}
\frac{b_{j}\cdots b_{N-1}}{a_{j}\cdots a_{N-1}}z_{j}+q^{-Nx}z_{0}=0, 
\quad (\gamma=1).
\end{array}
\right.
\label{contracond}
\ena
However, it can be checked that $Q_{n}$ does not satisfy (\ref{contracond}) 
if $n>2$ and $f_{n}\not=0$. 
This contradicts (\ref{mustrel}).

Hence $n=0$ or $1$, and $\{\Psi_{1}, \Psi_{2}\}$ is a basis.
\qed
 
\section{The $q$-hypergeometric function at $q^{N}=1$}
We recall the definition of the $q$-hypergeometric function at $|q|=1$. 

Set $q=e^{2\pi i\omega}, \omega >0$. 
The $q$-hypergeometric function of the Barnes type is given as follows \cite{NU, T}:
\bea
\Psi(\alpha, \beta, \gamma; x):=
\frac{\langle \alpha \rangle \langle \beta \rangle}{\langle 1 \rangle \langle \gamma \rangle}
\int_{C}q^{xz}
\frac{\langle z+1+\frac{1}{\omega} \rangle \langle z+\gamma \rangle}
{\langle z+\alpha \rangle \langle z+\beta \rangle}dz.
\label{defhyp}
\ena
Here the function $\langle z \rangle$ is defined by
\bea
\langle z \rangle :=
\exp{(\frac{\pi i}{2}\left( (1+\omega)z-\omega z^{2} \right))}
S_{2}(z | 1, \frac{1}{\omega}),
\ena
where $S_{2}(z)$ is the double sine function. 
We refer the reader to \cite{JM} for the double sine function.
The contour $C$ is the imaginary axis $(-i \infty, i\infty)$ 
except that the poles at
\bea
-\alpha +{\Bbb Z}_{\leq 0}+\frac{1}{\omega}{\Bbb Z}_{\leq 0}, \quad
-\beta +{\Bbb Z}_{\leq 0}+\frac{1}{\omega}{\Bbb Z}_{\leq 0} 
\ena
are on the left of $C$ and the poles at
\bea
{\Bbb Z}_{\geq 0}+\frac{1}{\omega}{\Bbb Z}_{\geq 0}, \quad 
-\gamma +{\Bbb Z}_{\geq 1}+\frac{1}{\omega}{\Bbb Z}_{\geq 1}
\ena
are on the right of $C$.
The integral (\ref{defhyp}) is absolutely convergent if
\bea
0< {\rm Re} x < 1+\frac{1}{\omega}+{\rm Re} \gamma -{\rm Re} \alpha -{\rm Re} \beta.
\label{conv}
\ena
Then the function $\Psi(x)$ satisfies the equation (\ref{qhyp}) 
with $q=e^{2\pi i\omega}$.

Now we consider $\Psi(x)$ under the condition (\ref{condition}).
We also assume that 
\bea
\alpha+\beta \le N-\gamma.
\label{condition2}
\ena
Then the integral (\ref{defhyp}) converges if $0<{\rm Re}x<1$ 
because $\omega=1/N$ and (\ref{conv}). 

\begin{thm} \label{th3}
Under the conditions (\ref{condition}) and (\ref{condition2}), 
the $q$-hypergeometric function $\Psi$ satisfies $\Psi \in {\cal S}$ and
is represented explicitly as follows:

1. $\gamma \le \beta \le \alpha$ case.
\bea
\Psi=\frac{1}{1-q^{Nx}}\Psi_{1}.
\label{rel1}
\ena

2. $\beta < \gamma \le \alpha$ case.
\bea
\Psi=\frac{1}{1-q^{Nx}} \left\{
\Psi_{1}+
\frac{(q)_{\gamma-1}}{(q)_{\alpha-1}(q)_{\beta-1}}
\frac{(q)_{\alpha-\gamma}(q)_{N-\gamma+\beta}}{(q)_{N-\gamma+1}}
\Psi_{2} \right\}.
\label{rel2}
\ena

3. $\beta \le \alpha < \gamma$ case.
\bea
&&
\Psi=\frac{1}{1-q^{Nx}}
\frac{(q)_{\gamma-1}}{(q)_{\alpha-1}(q)_{\beta-1}}
\frac{1}{N}
\frac{(q)_{N-\gamma+\alpha}(q)_{N-\gamma+\beta}}{(q)_{N-\gamma+1}}
\left( C_{\alpha, \beta}^{\gamma}\Psi_{1}-\Psi_{2} \right),  \label{rel3} \\
&&
C_{\alpha, \beta}^{\gamma}:=
1-(\sum_{j=1}^{N-\gamma+1}+\sum_{j=N-\gamma+\alpha+1}^{N-1}+\sum_{j=N-\gamma+\beta+1}^{N-1})
\frac{q^{j}}{1-q^{j}}. \label{defC}
\ena
Here $\left\{\Psi_{1}, \Psi_{2}\right\}$ in (\ref{rel1}), (\ref{rel2}) and (\ref{rel3}) 
is the basis of the space of solutions given in 
(\ref{base1}), (\ref{base2}) and (\ref{base3}), respectively.
\end{thm}

\proof
Let us calculate the integral in (\ref{defhyp}).
We denote by $\Phi(z)$ the integrand of (\ref{defhyp}):
\bea
\Phi(z):=q^{xz}
\frac{\langle z+1+\frac{1}{\omega} \rangle \langle z+\gamma \rangle}
{\langle z+\alpha \rangle \langle z+\beta \rangle}.
\ena
By using
\bea
\frac{\langle z+N \rangle}{\langle z \rangle}=
\frac{1}{1-e^{2\pi i z}},
\ena
we see the following under the condition (\ref{condition}):
\bea
\Phi(z+N)=q^{Nx}\Phi(z).
\ena
Hence, we have
\bea
(1-q^{Nx})\int_{C}\Phi(z)dz=
\left( \int_{C}-\int_{C+N} \right)\Phi(z)dz
\label{shift}
\ena
{}From (\ref{condition}), 
we can take the line $-\frac{1}{2}+i{\Bbb R}$ as the contour $C$. 
Then the right hand side of (\ref{shift}) is given by the sum 
of residues at $z=0, \cdots , N-1$.
Therefore, we get
\bea
\Psi(x)=
\frac{-2\pi i}{1-q^{Nx}}
\frac{\langle \alpha \rangle \langle \beta \rangle}
{\langle 1 \rangle \langle \gamma \rangle}
\sum_{k=0}^{N-1}{\rm res}_{z=k}\Phi(z).
\ena

By using
\bea
\frac{\langle z+1 \rangle}{\langle z \rangle}=
\frac{1}{1-q^{z}},
\ena
we can represent the function $\Phi(z)$ in terms of $q^{z}$ and
calculate residues of this function explicitly.

It is easy to see that all the poles at $z=0, \cdots N-1$ are simple 
if $\beta \le \alpha$ and $\gamma \le \alpha$. 
Then we find the formulae (\ref{rel1}) and (\ref{rel2}) easily
by using
\bea
(1-q)(1-q^{2})\cdots (1-q^{N-1})=N.
\label{basicrel}
\ena

Let us consider the case of $\beta \le \alpha < \gamma$. 
The poles at $z=0, \cdots , N-\gamma$ and $z=N-\alpha+1, \cdots , \beta$ are simple, 
and it is easy to calculate residues at these poles. 
The result is as follows:
\bea
{\rm res}_{z=k}=
-\frac{1}{N}\frac{(q)_{N-\gamma+\alpha}(q)_{N-\gamma+\beta}}{(q)_{N-\gamma+1}}
(1-q^{N-\gamma+1})
\frac{a_{k+1}\cdots a_{N-\gamma}}{b_{k}\cdots b_{N-\gamma}}, 
\quad (k=0, \cdots , N-\gamma)
\ena
and
\bea
{\rm res}_{z=k}&=&
-\frac{1}{N}\frac{(q)_{N-\gamma+\alpha}(q)_{N-\gamma+\beta}}{(q)_{N-\gamma+1}} \no \\
&\times&
(1-q^{\beta-\alpha})
\frac{b_{N-\gamma-1}\cdots b_{N-\alpha-1}}{a_{N-\gamma+2}\cdots a_{N-\alpha+1}}
\frac{b_{N-\alpha+1}\cdots b_{k-1}}{a_{N-\alpha+2}\cdots a_{k}}, 
\quad (k=N-\alpha+1, \cdots , \beta).
\ena
Here we used (\ref{basicrel}).
 
Next we calculate residues at $z=N-\gamma+1, \cdots , N-\alpha$. 
Note that these points are double poles. 
For $k=N-\gamma+1, \cdots , N-\alpha$, 
we have
\bea
&&
{\rm res}_{z=k}\Phi(z)={} \no \\
&&
=q^{kx}{\rm res}_{z=0}
\left(
\frac{q^{zx}}{(1-q^{z})^{2}}
\prod_{j=1}^{k}\frac{1}{1-q^{z+j}}
\prod_{j=\alpha+k}^{N-1}\frac{1}{1-q^{z+j}} \!\!\!\!\!
\prod_{j=1}^{k+\gamma-1-N}\!\!\!\!\! \frac{1}{1-q^{z+j}}
\prod_{j=\beta+k}^{N-1}\frac{1}{1-q^{z+j}}\right).
\label{res}
\ena
By substituting
\bea
&&
q^{zx}=1+\frac{2\pi i}{N}xz+o(z), \no \\
&&
\frac{1}{1-q^{z+j}}=
\left\{
\begin{array}{l}
\displaystyle 
-\frac{N}{2\pi i}z^{-1}+\frac{1}{2}+o(1), 
\quad (j\equiv 0  \, {\rm mod}\, N), \\
\displaystyle
\frac{1}{1-q^{j}}+\frac{2\pi i}{N}\frac{q^{j}}{(1-q^{j})^{2}}z+o(z), 
\quad (j \not\equiv 0 \, {\rm mod}\, N), 
\end{array}
\quad {\rm as} \quad z \to 0,
\right.
\ena
we find
\bea
&&
(\ref{res})=q^{kx}{\rm res}_{z=0}
\left( (\frac{N}{2\pi i})^{2}D_{k}z^{-2}+
  \frac{N}{2\pi i}D_{k}(x-E_{k})z^{-1}+o(1) \right) \no \\
&&
=\frac{N}{2\pi i}D_{k}(x-E_{k})q^{kx},
\ena
where
\bea
&&
D_{k}=
\prod_{j=1}^{k}\frac{1}{1-q^{j}}
\prod_{j=\alpha+k}^{N-1}\frac{1}{1-q^{j}}
\prod_{j=1}^{k+\gamma-1-N}\frac{1}{1-q^{j}}
\prod_{j=\beta+k}^{N-1}\frac{1}{1-q^{j}}, \no \\
&&
E_{k}=1-
\left(
\sum_{j=1}^{k}+\sum_{j=\alpha+k}^{N-1}+\sum_{j=1}^{k+\gamma-1-N}+\sum_{j=\beta+k}^{N-1} 
\right)
\frac{q^{j}}{1-q^{j}}.
\ena
By using (\ref{basicrel}), we see
\bea
D_{k}=\frac{1}{N^{2}}
\frac{(q)_{N-\gamma+\alpha}(q)_{N-\gamma+\beta}}{(q)_{N-\gamma+1}}
\frac{b_{N-\gamma+1}\cdots b_{k-1}}{a_{N-\gamma+2}\cdots a_{k}}.
\ena
Moreover, we have
\bea
&&
E_{N-\gamma+1}=C_{\alpha, \beta}^{\gamma},  \\
&&
E_{k}-E_{k-1}=\frac{1-q^{\gamma+2k-1}}{a_{k}}-\frac{1-q^{\alpha+\beta+2(k-1)}}{b_{k-1}},
\ena
where $C_{\alpha, \beta}^{\gamma}$ is given by (\ref{defC}).
Hence we find
\bea
E_{k}=C_{\alpha, \beta}^{\gamma}+
\sum_{j=N-\gamma+2}^{k}\left(
\frac{1-q^{\gamma+2j-1}}{a_{j}}-\frac{1-q^{\alpha+\beta+2(j-1)}}{b_{j-1}}
\right).
\ena

{}From this calculation, we get the relation (\ref{rel3}).
\qed \newline

Theorem \ref{th3} implies that 
the $q$-hypergeometric function of the Barnes type has 
logarithmic singularity in the case that $q^{N}=1$ and $\beta \le \alpha < \gamma$.

\end{document}